\documentclass[preprint, 11pt]{elsarticle}

\usepackage[T1]{fontenc}
\usepackage[utf8]{inputenc}
\usepackage[hmargin=0.9in]{geometry}
\usepackage[babel,german=quotes]{csquotes}
\usepackage{subcaption}
\usepackage{graphicx}
\usepackage[colorlinks=true,citecolor=blue,urlcolor=blue]{hyperref}
\usepackage{caption}
\usepackage{amsmath}
\usepackage{amsthm}
\usepackage{amssymb}
\usepackage{bm, amsfonts}
\usepackage{xcolor,soul,framed}
\usepackage{fixmath}
\usepackage{tikz}
\usetikzlibrary{patterns, arrows.meta}
\usepackage{bm}
\usepackage{makecell}
\usepackage{titlesec}
\usepackage{multirow}
\usepackage{indentfirst}
\usepackage{mathtools}
\usepackage{float}
\usepackage[normalem]{ulem}
\useunder{\uline}{\ul}{}
\usepackage{placeins}

\usepackage{booktabs} 
\usepackage{siunitx}  
\usepackage{adjustbox} 
\usepackage{tikz}

\newcommand{\R}{\mathbb{R}}
\newcommand{\N}{\mathbb{N}}

\colorlet{shadecolor}{green}

\titleformat{\paragraph}
{\normalfont\normalsize\itshape}{\theparagraph}{1em}{}
\titlespacing*{\paragraph}
{0pt}{3.25ex plus 1ex minus .2ex}{1.5ex plus .2ex}

\newdefinition{remark}{Remark}

\makeatletter
\def\ps@pprintTitle{%
\let\@oddhead\@empty
\let\@evenhead\@empty
\def\@oddfoot{
\footnotesize\itshape
\ifx\@journal\@empty Elsevier
\else\@journal\fi
\hfill\today
}%
\let\@evenfoot\@oddfoot}
\makeatother

\newcommand{\D}{\mathrm{d}}
\newcommand{\td}[2]{\frac{\D #1}{\D #2}}

\newcommand{\pd}[2]{\frac{\partial #1}{\partial #2}}
\newcommand{\dx}{\mathrm{d}\mathbf{x}}
\newcommand{\ds}{\mathrm{d}\mathbf{s}}
\newcommand{\bu}{\mathbf{u}}
\newcommand{\bU}{\mathbf{U}}
\newcommand{\bomega}{\boldsymbol{\omega}}
\newcommand{\bx}{\mathbf{x}}
\newcommand{\bX}{\mathbf{X}}

\newcommand{\bT}{\mathbf{T}}
\newcommand{\hbu}{\hat{\mathbf{u}}}
\newcommand{\hp}{\hat p}

\newcommand{\beq}{\begin{equation}}
\newcommand{\eeq}{\end{equation}}

\titleformat{\subsection}{\bfseries\small}{\thesubsection}{1em}{}

\begin{document}

\begin{frontmatter}
  \title{A Chimera domain decomposition method with weak Dirichlet-Robin coupling 
    for finite element simulation of particulate flows}
\author{Raphael Münster, Otto Mierka,  Dmitri Kuzmin$^*$, Stefan Turek}
\ead{\{raphael.muenster,otto.mierka,dmitri.kuzmin,stefan.turek\}@math.tu-dortmund.de}
\cortext[cor1]{Corresponding author}

\address{Institute of Applied Mathematics (LS III), TU Dortmund University\\ Vogelpothsweg 87, D-44227 Dortmund, Germany}

\journal{Mathematics and Computers in Simulation}

\begin{abstract}   
    We introduce a new multimesh finite element method for direct numerical simulation of incompressible particulate flows. The proposed approach falls into the category of overlapping domain decomposition / Chimera / overset grid meshes. In addition to calculating the velocity and pressure of the fictitious fluid on a fixed background mesh, we solve the incompressible Navier-Stokes equations on body-fitted submeshes that are attached to moving particles. The submesh velocity and pressure are used to calculate the hydrodynamic forces and torques acting on the particles. The coupling with the background velocity and pressure is enforced via (i) Robin-type boundary conditions for an Arbitrary-Lagrangian-Eulerian (ALE) formulation of the submesh problems and (ii) a Dirichlet-type distributed interior penalty term in the weak form of the background mesh problem. The implementation of the weak Dirichlet-Robin coupling is discussed in the context of discrete projection methods and finite element discretizations. Detailed numerical studies are performed for standard test problems involving fixed and moving immersed objects. A comparison of Chimera results with those produced by fictitious boundary methods illustrates significant gains in the accuracy of drag and lift approximations.
\end{abstract}


\end{frontmatter}


\bigskip

\section{Introduction}
\label{sec:intro}

Numerical methods for direct numerical simulation (DNS) of incompressible flows around
moving rigid particles can be classified into fixed mesh and deforming mesh
approaches \cite{review2012}. Prominent representatives of the latter family
include immersed boundary methods \cite{lee2011,peskin,uhlmann} and 
fictitious domain formulations in which the rigid body motion inside
the particles is enforced using distributed Lagrange multipliers
(DLM, \cite{glowinski1,glowinski2,patankar2000}), fictitious (surrogate, shifted,
unfitted) boundary methods \cite{sbm2018,vonwahl-diss,WanTurek2006a,WanTurek2006b,WanTurek2007a}, subspace projections \cite{steffen,rod}, and other non-DLM alternatives
\cite{minev2007,sharma2005}. The resolving power of such approaches can
be greatly improved by using mesh deformation techniques
\cite{steffen,fbm2012,WanTurek2007b} or overlapping domain decomposition
with moving submeshes
\cite{multimesh2020,henshaw1994,henshaw2001,codina,multimesh2019}. The
coupling conditions in the overlap region can again be enforced in
different ways. Houzeaux and Codina \cite{codina} propose a Chimera
method with Dirichlet/Neumann(Robin) coupling such that Dirichlet-type
conditions are strongly enforced at nodal points inside and around
immersed objects. The multimesh finite element method developed by
Dokken et al. \cite{multimesh2020} achieves the coupling effect by
incorporating suitable stabilization terms into the discontinuous
Galerkin (DG) weak forms of interacting subproblems.

A common disadvantage of all fixed-mesh algorithms that constrain the
velocity of a fictitious fluid at discrete locations is the lack of
continuous dependence on the data. Indeed, a small displacement of
a particle can activate or deactivate the Dirichlet constraint.
Moreover, the volume of the constrained region changes abruptly
leading to nonphysical temporal oscillations in the pressure field and,
as a consequence, in forces acting on the particles. In the
present paper, we cure this deficiency of Chimera-type domain
decomposition methods by using a distributed interior penalty term
instead of strongly imposed nodal constraints. The weak form of
our background mesh problem differs from the unfitted finite
element method presented in \cite{multimesh2020} in the structure of
employed stabilization terms and in the way in which they are incorporated
into the discrete projection method for the Navier--Stokes system.
Our algorithm is simpler than the DG
method from \cite{multimesh2020} and requires fewer degrees of freedom because
we are using a continuous approximation to the velocity field. The
results of our numerical experiments show that our modification of
the Chimera method from \cite{codina} 
is robust and capable of delivering accurate DNS results
at a fraction of the cost that a fixed-mesh fictitious boundary
method would require.

The rest of this paper is organized as follows. A fictitious
domain formulation of the particulate flow problem that we consider
is presented in Section \ref{sec:fdm}. The Chimera domain decomposition method
introduced in Section  \ref{sec:chimera} leads to a system of weakly coupled
Navier-Stokes problems. In Section~\ref{sec:fem}, we incorporate Dirichlet
constraints into the weak form of the background mesh problem
using the distributed interior penalty approach. In Section \ref{sec:proj},
we show how we update the velocity and pressure using the framework
of discrete projection methods. The new Chimera penalty method for solving
the coupled problem is summarized in Section \ref{sec:summary} and
compared with a modification in which the Dirichlet constraints are
imposed strongly (as in the original method \cite{codina}).
The numerical behavior of both
versions is studied in Section \ref{sec:num}, and some conclusions
are drawn in Section \ref{sec:concl}.

\section{Fictitious domain formulation}
\label{sec:fdm}

Let $\Omega\subset\R^d,\ d\in\{2,3\}$ be a fixed fictitious domain. In
our particulate flow model, $\Omega$ is filled
with an incompressible Newtonian fluid that carries
a suspension of $N_p$ rigid
particles (balls) $$B_k(t)=\{\mathbf{x}\in\R^d\,:\,|\mathbf{x}-\mathbf X_k(t)|< R_k\}.$$
We subdivide $\Omega$ into 
$\Omega_p(t):=\bigcup_{k=1}^{N_p}B_k(t)$
and the subdomain $\Omega_f(t):=\Omega\backslash\bar\Omega_p(t)$
occupied by the fluid.

The evolution of the fluid-particle
mixture is governed by the system \cite{fbm2012,WanTurek2006b}
\begin{subequations}\label{nse}
\begin{align}\label{nse1}
\rho_f\left(\pd{\mathbf{u}}{t}+\mathbf u\cdot\nabla\mathbf u\right)
&= -\nabla p +
\nabla\cdot (2\mu_f\mathbf{D}(\mathbf u))\quad \mbox{in}\ \Omega_f(t),\\
  \nabla\cdot \mathbf u&=0\quad \mbox{in}\
  \Omega,\label{nse2}\\
  \mathbf{u}&=\bU\quad \mbox{on}\ \bar\Omega_p(t)\label{nse3}
\end{align}
\end{subequations}
of generalized incompressible Navier--Stokes equations. Here $\bu$ is
the velocity of the fluid, $\bU$ is the velocity of rigid body
motion inside the particles, $p$ is the pressure, and $$\mathbf{D}(\mathbf u)
=\frac12(\nabla\bu+\nabla\bu^\top)$$ is the deformation rate tensor. The constant
density and dynamic viscosity of the fluid phase are denoted by
$\rho_f$ and $\mu_f$, respectively. For $\bx\in\bar B_k(t)$, the
velocity of the fictitious fluid is given by
\beq\label{rbm}
\bU(\bx,t)=\bU_k(t)+\bomega_k(t)\times(\bx-\bX_k),
\eeq
where $\bU_k$ is the translational velocity and $\bomega_k$ is the
angular velocity of the $k$th particle. The density, volume,
and moment of inertia tensor of this particle are denoted by
$\rho_k,\ V_k$, and $\mathbb{I}_k$, respectively.

For simplicity, we assume that the
particles do not collide with each other or with solid
walls. Therefore, we do not include repulsive or lubrication
forces in the Newton--Euler equations
\begin{subequations}\label{newton}
\begin{align}
  \rho_pV_k\td{\bU_k}{t}&=\mathbf{F}_k
  +(\rho_p-\rho_f)V_k\mathbf{g},\\
  \mathbb{I}_k\td{\bomega_k}{t}&=
  \bT_k-\bomega_k\times(\mathbb{I}_k\bomega_k),
\end{align}
\end{subequations}
where $\mathbf{g}$ is the gravitational acceleration.
 The hydrodynamic force $\mathbf{F}_k$
 and torque $\mathbf{T}_k$ are defined by
 \begin{subequations}
   \label{fktkdef}
\begin{align}   
\mathbf{F}_k&=-\int_{\partial B_k}\boldsymbol{\sigma}
\mathbf{n}\ds,\\
\bT_k&=-\int_{\partial B_k}(\bx-\bX_k)\times
(\boldsymbol{\sigma}\mathbf{n})\ds.
\end{align}
\end{subequations}
Here $\mathbf{n}$ denotes the unit outward normal
and $\boldsymbol{\sigma}=
-p\mathbf{I}+2\mu_f\mathbf{D}(\mathbf u)$, where 
$\mathbf{I}$ is the identity tensor.

The problem statement is completed by imposing appropriate
initial and boundary conditions. By default, we prescribe
a given velocity profile at the inlet, the zero-stress
condition at the outlet, and the no-slip condition on
solid walls. The rigid body motion constraint \eqref{nse3}
defines both the no-slip Dirichlet boundary data
for $\partial\Omega_p(t)$ and the velocity of the
fictitious fluid contained in $\Omega_p(t)$.

\section{Chimera domain decomposition}
\label{sec:chimera}

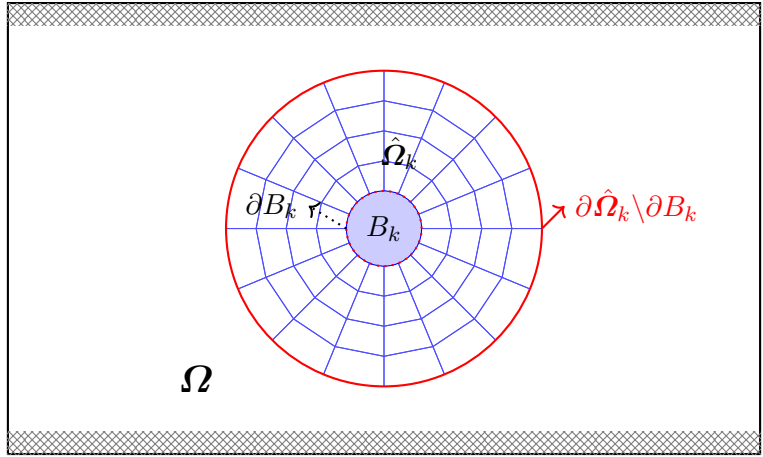
\begin{figure}[t]
\begin{center}
\begin{tikzpicture}

\def\rstep{0.4}
\def\rmax{2.0}
\def\divisions{16}
\def\rinner{0.5}

\draw[thick] (-5,-3) rectangle (5,3);

\fill[pattern=crosshatch, pattern color=black!50] (-5,-3) rectangle (5,-2.7);

\fill[pattern=crosshatch, pattern color=black!50] (-5,2.7) rectangle (5,3);

\node at (-2.5,-2) {\Large$\boldsymbol{\Omega}$};

\foreach \r in {0.5,0.9,1.3,1.7} {
  \pgfmathsetmacro{\rnext}{\r + \rstep}
  \foreach \i in {0,...,15} {
    \pgfmathsetmacro{\anglea}{\i * 360 / \divisions}
    \pgfmathsetmacro{\angleb}{(\i+1) * 360 / \divisions}

    \draw[blue!70]
      (\anglea:\r) -- (\angleb:\r) -- (\angleb:\rnext);

    \ifdim\rnext pt<2.0pt
      \draw[blue!70] (\angleb:\rnext) -- (\anglea:\rnext);
    \fi

    \draw[blue!70] (\anglea:\rnext) -- (\anglea:\r);
  }
}

\filldraw[blue!20, draw=blue] (0,0) circle (\rinner);

\draw[red, dotted, thick] (0,0) circle (\rinner);

\draw[red, thick] (0,0) circle (2.1);

\node at (0,0) {$B_k$};
\node at (0.2,1.0) {$\hat{\boldsymbol{\Omega}}_k$};

\draw[->,thick,red] (\rmax + 0.1,0) -- (\rmax+0.4,0.3)
  node[right, red] {$\partial \hat{\boldsymbol{\Omega}}_k \backslash \partial B_k$};

\draw[->, thick, dotted] (-\rinner,0) -- (-\rinner-0.5,0.3)
  node[left] {$\partial B_k$};

\end{tikzpicture}
\end{center}
    \caption{Subdomains of the Chimera domain decomposition method.}
    \label{fig:chimera_sketch}
\end{figure}

In our multimesh numerical method for solving the coupled problems
\eqref{nse} and \eqref{newton}, we calculate
$\mathbf{F}_k$ and $\bT_k$ using finite element approximations
to $\bu$ and $p$ on the body-fitted subdomains $$\hat\Omega_k(t)
=\{\mathbf{x}\in\R^d\,:\, R_k<|\mathbf{x}-\mathbf X_k(t)|< R_k+H_k\},$$
which are embedded into $\Omega$ as shown in Fig.~\ref{fig:chimera_sketch}. In this work, we assume
that $\hat\Omega_k(t)\cap B_j(t)=\emptyset$ for $k\ne j$. The
general case of overlapping domains is considered in
\cite{multimesh2020,multimesh2019}.

If the particle $B_k(t)$ were a planet,
the associated subdomain $\hat\Omega_k(t)$ could be interpreted
as the atmosphere of that planet. We denote the atmospheric
velocity  and pressure fields by $\hbu$ and $\hp$, respectively.
The parameter $H_k>0$ determines the width of the
atmospheric layer around $B_k(t)$.

Following Houzeaux and Codina \cite{codina}, we perform \emph{iteration
by subdomains} using coupling conditions of Dirichlet--Robin type.
That is, the background fields $\bu$ and $p$ influence the
solution of
\begin{subequations}\label{snse}
\begin{align}\label{snse1}
\rho_f\left(\pd{\hbu}{t}+\hbu\cdot\nabla\hbu\right)
&= -\nabla \hp +
\nabla\cdot (2\mu_f\mathbf{D}(\hbu))\quad \mbox{in}\ \hat\Omega_k(t),\\
  \nabla\cdot \hbu&=0\quad \mbox{in}\
  \hat\Omega_k(t),\label{snse2}\\
  \hbu&=\bU\quad \mbox{on}\ \partial B_k(t),\label{snse3} \\
  \hat{\boldsymbol{\sigma}}\mathbf{n}-\alpha(\hbu\cdot\mathbf{n})\hbu&=
  \boldsymbol{\sigma}\mathbf{n}-\alpha(\bu\cdot\mathbf{n})\bu
  \quad \mbox{on}\ \partial\hat\Omega_k(t)\backslash
  \partial B_k(t)\label{snse4}
\end{align}
\end{subequations}
by providing the data of the Robin/Neumann boundary condition \eqref{snse4},
where $\alpha\ge 0$ is an interior penalty parameter and
$\hat{\boldsymbol{\sigma}}=-\hat p\mathbf{I}+2\mu_f\mathbf{D}(\hbu)$
is the total stress. The consistency relation
\beq\label{dircond}
\bu=\hbu\quad\mbox{on}\ \bar{\hat \Omega}_k(t)
\eeq
is satisfied by exact solutions to \eqref{nse} and \eqref{snse}.
Hence, it is appropriate to constrain a numerical approximation
to the background velocity
$\bu$ using a strong or weak form of the Dirichlet condition
\eqref{dircond}.

\section{Finite element discretization}
\label{sec:fem}

We discretize the velocity and pressure
fields in space
using the inf-sup stable $\mathbb{Q}_2$-$\mathbb{P}_1^{\rm disc}$
finite element pair. The atmospheric subproblems \eqref{snse}
are solved on body-fitted meshes using an Arbitrary-Lagrangian-Eulerian
(ALE) formulation in the reference frame moving with the constant
mesh velocity $\mathbf w_k(\bx,t)=\bU_k(t)$. Details
of such moving mesh finite element (FE) methods can be found elsewhere
\cite{steffen,maury1999,WanTurek2007b} and are not discussed here. The
ALE-FE approximation to $(\hbu,\hp)$ is denoted by $(\hbu_h,\hp_h)$.

The fictitious domain problem \eqref{nse} is solved using a fixed
background mesh $\mathcal T_h$ that consists of quadrilaterals (in 2D)
or hexahedra (in 3D). Instead of enforcing the Dirichlet-type velocity
constraints \eqref{nse3} and \eqref{dircond} strongly at discrete locations
(as in \cite{codina,rod,WanTurek2006a}), we incorporate them into a
weak form of problem \eqref{nse} using (a discrete counterpart of) the
distributed interior penalty term
\beq\label{intpen}
s(\hbu,\bU;\bu,\mathbf v)=\gamma_{\max}\sum_{k=1}^{N_p}\left[\int_{\hat\Omega_k(t)}
  \beta_k(\bu-\hbu)\cdot \mathbf v\dx+
  \int_{B_k(t)}(\bu-\bU)\cdot \mathbf v\dx
  \right]
\eeq
that depends on a penalty parameter $\gamma_{\max}\gg 1$ and a damping
function $\beta_k:\hat\Omega_k\times\R_0^+\to [0,1]$ such as\footnote{
The nonnegative function $\beta_k$ determines the strength of velocity penalization.
It should vanish in a neighborhood of the interface
$\Gamma_k(t)=\partial\hat\Omega_k(t)\backslash \partial B_k(t)$
to avoid interference with the Robin boundary condition 
\eqref{snse4}.}
$$\beta_k(\bx,t)=\min\left(1,\max\left(0,
\frac{R_k+0.75H_k-|\bx-\bX_k(t)|}{0.25H_k}\right)\right).
$$
The velocity $\bU$ of rigid body motion
is defined by \eqref{rbm}. In the discrete
version of \eqref{intpen},
integration is performed over the regions $\hat \Omega_{k,h}$
and $B_{k,h}$ that are covered by / enclosed by the ALE submesh.

\begin{figure}[h]
    \centering
    \includegraphics[width=0.85\textwidth,trim=50 490 0 0,clip]{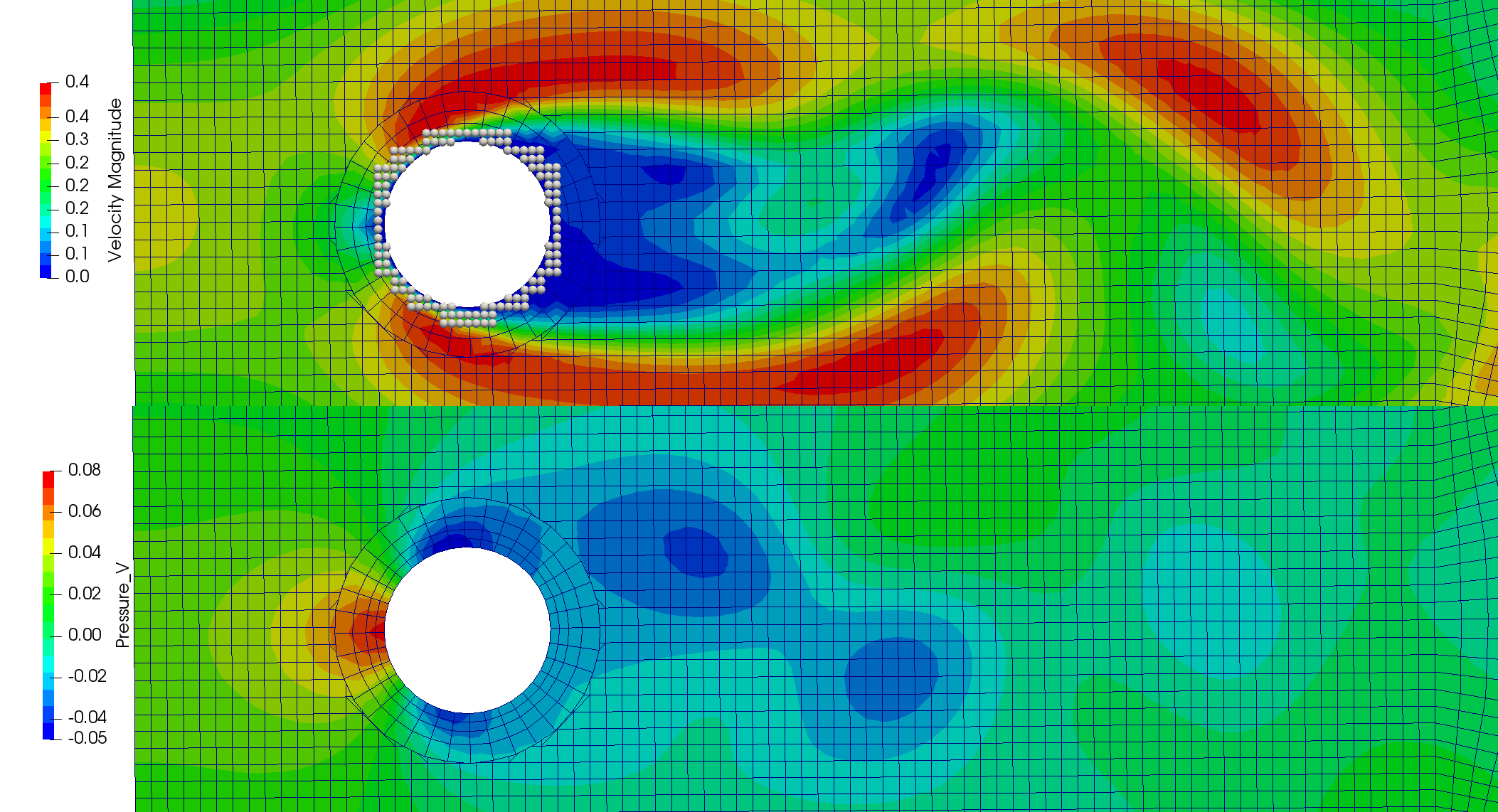} 
    \vskip0.75cm
    
    \includegraphics[width=0.85\textwidth,trim=50 0 0 570,clip]{fL4L3.0010.png}
    \vskip0.25cm
    
    \caption{Chimera-type domain decomposition for a flow around a cylinder
      configuration.
      Velocity (top) and pressure (bottom) distribution on a uniform
      background mesh and an overlapping body-fitted submesh.}

    \label{fig:Chimera_snapshot}
\end{figure}

Let $\mathbf{V}_h$ and $Q_h$ denote the finite
element spaces for the velocity and pressure approximations
on $\bar\Omega$, respectively. We
 seek $\bu_h\in \mathbf{V}_h$ and $p_h\in Q_h$ such that
\begin{subequations}\label{wf}
  \begin{align}\nonumber
    \int_{\Omega}
    \rho_f\left(\pd{\mathbf{u}_h}{t}+\mathbf u_h\cdot\nabla\mathbf u_h\right)
    \cdot\mathbf v_h\dx&+
    \frac{\mu_f}2\int_\Omega \mathbf D(\bu_h):\mathbf D(\mathbf v_h)\dx
    +s_h(\hbu_h,\bU_h;\bu_h,\mathbf v_h)\label{wf1}\\
    &-\int_\Omega p_h\nabla\cdot\mathbf v_h\dx=0, \qquad  \mathbf v_h\in
    \mathbf{V}_h,\\
    &\phantom{-}\int_{\Omega}q_h \nabla\cdot \mathbf u_h\dx =0,
    \qquad \mathbf q_h\in Q_h.\label{wf2}
\end{align}
\end{subequations}
Figure \ref{fig:Chimera_snapshot} shows 
a snapshot of the region in which the ALE
mesh attached to $B_k(t)$
intersects the fixed background mesh. The former provides the
Dirichlet data $\hbu_h$ for calculating $s_h(\hbu_h,\bU_h;\bu_h,\mathbf v_h)$
defined by \eqref{intpen}. The latter provides the Robin boundary data
$\boldsymbol{\sigma}_h\mathbf{n}-\alpha(\bu_h\cdot\mathbf{n})\bu_h$
for updating $\hbu_h$. The hydrodynamic force $\mathbf F_k$ and torque
$\mathbf T_k$ are calculated using 
$\hat{\boldsymbol{\sigma}}_h=-\hat p_h\mathbf{I}+2\mu_f\mathbf{D}(\hbu_h)$
to approximate $\boldsymbol{\sigma}$ in \eqref{fktkdef}. The
high accuracy of the submesh approximation to the traction
$\hat{\boldsymbol{\sigma}}_h\mathbf n$
on $\partial B_k$ is a key advantage of the Chimera approach
compared to one-mesh fictitious
domain methods \cite{rod,minev2007,WanTurek2006a}.

\begin{remark}
A remarkable property of the weak form \eqref{wf} is that the
interior penalty term \eqref{intpen} is well defined even for
overlapping `atmospheres' $\hat\Omega_k(t)$ and  $\hat\Omega_j(t)$
of non-overlapping particles $B_k(t)$ and $B_j(t)$. The
use
of strongly imposed Dirichlet coupling conditions at background
mesh nodes belonging to
the overlap region $\hat\Omega_k(t)\cap\hat\Omega_j(t)$
would require  artificial averaging of
submesh data. The proposed approach performs such
averaging automatically by using penalty parameters depending on
 $\beta_k$ and $\beta_j$. An extension to the case 
$\hat\Omega_k(t)\cap B_j(t)\ne\emptyset$ is feasible but would
require adding a Dirichlet penalty term to the submesh
problem for $\hat\Omega_k(t)$. Such 
extensions can be
carried out following \cite{multimesh2020,multimesh2019}.
\end{remark}

\section{Fractional-step method}
\label{sec:proj}

Let $t^n=n\Delta t$, where $\Delta t$ is a constant time step and
$n\in\N_0$. Approximate solutions at the time level
$t^n$ will be referred to using the superscript $n$.
At the beginning of each time step, we
use the old submesh approximation
$(\hbu_h^{n},\hp_h^{n})$ to calculate the surface integrals
\eqref{fktkdef}. The old background mesh approximation
$(\bu_h^{n},p_h^{n})$ is used to calculate the data
of the Robin boundary condition
\eqref{snse4}. Then we update the positions of the
particles and solve the discrete saddle point problems
associated with ALE submeshes for $\hat\Omega_k(t)$.
These problems are small and can be solved efficiently, e.g.,
using the local Multilevel Pressure Schur Complement (MPSC) method
\cite{Turek1999,afc3,featflow}. The updated approximation
$(\hbu_h^{n+1},\hp_h^{n+1})$ can then be substituted into the interior
penalty term of problem \eqref{wf} for $(\bu_h^{n+1},p_h^{n+1})$.

Using the two-level $\theta$-scheme to discretize
\eqref{wf1} in time, we obtain a nonlinear system
of the form
\beq\label{saddle}
\begin{bmatrix}
  A(u^{n+1})+D^{n+1} & B\\
  B^T & 0
\end{bmatrix}
\begin{bmatrix}
u^{n+1}\\ p^{n+1}
\end{bmatrix}=
\begin{bmatrix}
f^n+g^{n+1}\\0
  \end{bmatrix}.
\eeq
The contribution of the interior penalty term
\eqref{intpen} is represented by
$D^{n+1}u^{n+1}-g^{n+1}$, where
$D^{n+1}$ is a symmetric positive semi-definite matrix. In
Appendix A below,
we explicitly define all matrices and vectors that
appear in \eqref{saddle}, as well as the lumped
approximation $M_L$ to the
consistent
mass matrix $M_C$ of the finite element
space $\mathbf{V}_h$. Below we use
$M_L$ for preconditioning purposes.
\medskip

We solve the discrete problem \eqref{saddle}
approximately using the following fractional-step
algorithm:
\begin{enumerate}
\item Solve the viscous Burgers system
  \beq\label{burgers}
            [A(\tilde u^{n+1})+D^{n+1}]
              \tilde u^{n+1}=f^n+g^{n+1}-Bp^n.
  \eeq
\item Solve the pressure Poisson system
  \beq\label{ppe}
  B^\top M_L^{-1}B(p^{n+1}-p^n)=\frac{1}{\Delta t}
  B^\top\tilde u^{n+1}.
  \eeq
\item Correct the intermediate velocity
  \beq\label{corr}
  [M_L+\Delta tD^{n+1}]u^{n+1}=[M_L+\Delta tD^{n+1}]
  \tilde u^{n+1}-\Delta tB(p^{n+1}-p^{n}).
   \eeq
\end{enumerate}
Note that we penalize the velocity not only in the
Burgers step but also in the final update. The
nonlinear system \eqref{burgers} is linearized about $u^n$ or
solved using a fixed-point
iteration method \cite{Turek1999}.

\begin{remark}
Rearranging \eqref{corr} and invoking \eqref{ppe}, we find that the
corrected velocity $u^{n+1}$ satisfies
\begin{align*}
  B^\top u^{n+1}&=\Delta tB^\top M_L^{-1}D^{n+1}(\tilde u^{n+1}-u^{n+1}).
\end{align*}
Hence, $u^{n+1}$ is approximately divergence-free. Moreover, the
residual of the constraint $B^\top u^{n+1}=0$
can be made as small as desired by performing additional outer
iterations (cf.~\cite{Turek1999,afc3}), in which the right-hand
side of \eqref{burgers}  is recalculated using the latest
approximation to $p^{n+1}$ instead of $p^n$.
\end{remark}

\section{Strong vs. weak Dirichlet coupling}
\label{sec:summary}

For comparison purposes, we implemented a Chimera multimesh method
in which the Dirichlet constraints are imposed strongly at \emph{hole nodes} in
$\hat B_k(t)$ and \emph{fringe nodes} in $\hat \Omega_k$. The terminology that
we use here is adopted from \cite{codina}.
A node $\mathbf x_i$ of a mesh cell crossed by $\partial B_k(t)$ is treated
as a hole node if $\mathbf x_i\in B_k(t)$ and as a fringe node if 
$\mathbf x_i\in\hat\Omega_k(t)$. Following
Houzeaux and Codina \cite{codina},
we prescribe the fictitious boundary
conditions $\mathbf u_h(\mathbf x_i)=\mathbf U(\mathbf x_i)$ and
$\mathbf u_h(\mathbf x_i)=\hat{\mathbf u}_h(\mathbf x_i)$ at hole
and fringe nodes, respectively. Details can be found in \cite{codina}.
For a better comparison of Chimera methods with weak and strong
Dirichlet coupling, the main algorithmic steps are summarized below
for both versions.
\pagebreak

The strong-form implementation will be referred to as \emph{Chimera-S}.
It represents a combination of the methods developed
in \cite{codina} and \cite{WanTurek2006a,WanTurek2006b,WanTurek2007a}.
The numerical solutions are updated as follows:
\begin{enumerate}
\item Solve the background mesh problem \eqref{nse} using
  fictitious Dirichlet boundary conditions at hole nodes
  in the first outer iteration and at hole+fringe nodes in subsequent updates.
\item Solve the submesh problems \eqref{snse} using interpolated background velocity $\bu_h$
  and pressure $p_h$ in the surface integral associated with the Robin boundary condition \eqref{snse4}.
\item Calculate $\mathbf{F}_k$ and $\mathbb{I}_k$ for each particle 
  using the submesh stress $\hat{\boldsymbol{\sigma}}_h$ in the
  surface integrals \eqref{fktkdef}.
\item Solve the discrete version of the Newton--Euler equations
  \eqref{newton}. Use the updated velocity $\bU_k$ and
  angular velocity $\bomega_k$ to determine the new position and
  orientation of particles.
\item Find the hole/fringe nodes corresponding to the updated
  position of particles, interpolate the submesh velocity $\hbu_h$ to
  these nodes of the background mesh 
  and go to Step 1 or exit.
\end{enumerate}
At least two outer iterations need to be performed
in this 
fictitious boundary Chimera method.
\medskip
  
The implementation labeled \emph{Chimera-W} incorporates the interior
penalty term \eqref{intpen} into
the discretized momentum equation
\eqref{wf1}. The proposed algorithm consists of the
following steps:
\begin{enumerate}
\item Calculate $\mathbf{F}_k$ and $\mathbb{I}_k$ for each particle 
  using the submesh stress $\hat{\boldsymbol{\sigma}}_h$ in the
  surface integrals \eqref{fktkdef}.
\item Solve the discrete version of the Newton--Euler equations
  \eqref{newton}. Use the updated velocity $\bU_k$ and
  angular velocity $\bomega_k$ to determine the new position and
  orientation of particles.
  \item Solve the submesh problems \eqref{snse} using interpolated background velocity $\bu_h$
  and pressure $p_h$ in the surface integral associated with the Robin boundary condition \eqref{snse4}.
\item Go to Step 1 if stronger coupling is desired.
\item Solve the background mesh problem \eqref{nse} 
  using the updated
  submesh velocity $\hbu_h$ in the interior penalty term
  $s_h(\hbu_h,\bU;\bu_h,\mathbf v_h)$ of
  the discretized momentum equation \eqref{wf1}.
\item Go to Step 1 if stronger coupling is desired.  
\end{enumerate}
In contrast to the Chimera-S version, this algorithm may be
terminated after one outer iteration.



\section{Numerical examples}
\label{sec:num}

In this section, we perform numerical studies of the Chimera-S
and Chimera-W domain decomposition methods presented in Section
\ref{sec:summary}. For comparison purposes, we also ran simulations
using the classical fictitious boundary method \cite{WanTurek2006a,WanTurek2006b}, in which no submeshes are used, and the background mesh is uniform. In this
version (FBM), the constraint of rigid body motion \eqref{rbm} is enforced strongly at the hole nodes inside the particles. Since we are interested in applications to particulate flows, some of our investigations are focused
on the accuracy of approximations to drag and lift forces.
The unfitted treatment of embedded boundaries in the FBM approach
results in large errors on coarse meshes. Moreover, hydrodynamic forces
typically exhibit oscillatory behavior. The objective of our
experiments is to show that the use of body-fitted submeshes
leads to significant improvements.

\subsection{DFG benchmark 2D-2}\label{FAC}
One of the most popular test problems in computational fluid dynamics (CFD) is
known as the DFG benchmark 2D-2 \cite{dfg}. It is designed to evaluate numerical algorithms for solving the incompressible Navier-Stokes equations under laminar flow conditions. Specifically, this benchmark provides reference data and computational setup for simulating the two-dimensional flow around a circular cylinder at the Reynolds number Re=100. The experimentally observed flow pattern exhibits periodic vortex shedding, known as the von K\'arm\'an vortex street. For an in-depth description of the benchmark, we refer the reader to \cite{dfg}. The most important quantities and parameters are defined below.

The computational domain is a two-dimensional rectangular channel with a circular obstacle representing a cross section of a cylinder. The length and height of the channel are 2.2 units and 0.41 units, respectively. The cylinder has a diameter of 0.1 units and is centered at (0.2, 0.2). 

Assuming that the laminar flow is fully developed at the inlet,
the parabolic velocity profile
\[
\mathbf{u}(0, y) = \left( \frac{4 U y (0.41 - y)}{0.41^2}, 0 \right)
\]
is prescribed on the inflow boundary ($x=0$). The maximum velocity $U$ is attained at $y = 0.205$. 

The drag force $F_d$ and the lift force $F_l$ exerted by the fluid on
the cylinder are given by
\[
\begin{pmatrix} F_d \\ F_l \end{pmatrix} =\frac{1}{\rho_f}
\int_{\partial B} \hat{\boldsymbol{\sigma}}_h\mathbf n
\ds,
\]
where $\hat{\boldsymbol{\sigma}}_h$ is the submesh approximation to
$\boldsymbol{\sigma}
=-p\mathbf{I}+2\mu_f\mathbf{D}(\bu)$,
and $\mathbf n$ is the unit outward normal to
the boundary $\partial B$ of the cylinder. Using $F_d$ and $F_l$, we compute the
drag and lift coefficients 
\begin{equation}
C_d = \frac{2 F_d}{U_{\text{mean}}^2 L}, \qquad C_l = \frac{2 F_l}{U_{\text{mean}}^2 L},
\label{eq:cd_cl}
\end{equation}
where $U_{\text{mean}}=1.0$ is the average velocity at the inlet
and $L=0.1$ is the characteristic length. Because of periodic
vortex shedding, the coefficients $C_d$ and $C_d$ vary over time.

\begin{figure}[h!]
    \centering
    \begin{subfigure}{0.8\textwidth}
        \centering
        \includegraphics[width=0.75\textwidth]{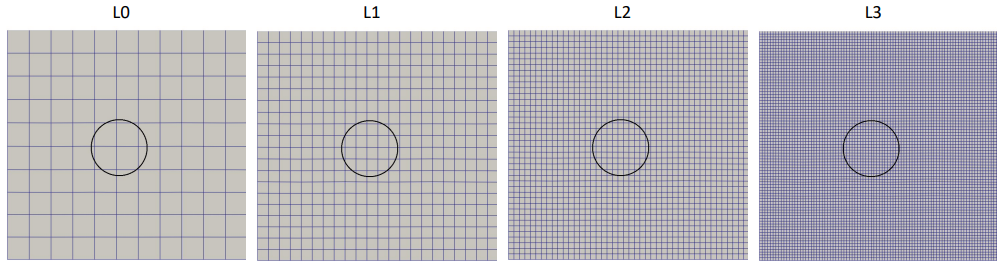} 
        \caption{Benchmark mesh with an immersed cylinder on different levels of resolution.}
        \label{fig:sub1}
    \end{subfigure}

    \vspace{1em} 

    \begin{subfigure}{0.8\textwidth}
        \centering
        \includegraphics[width=0.75\textwidth]{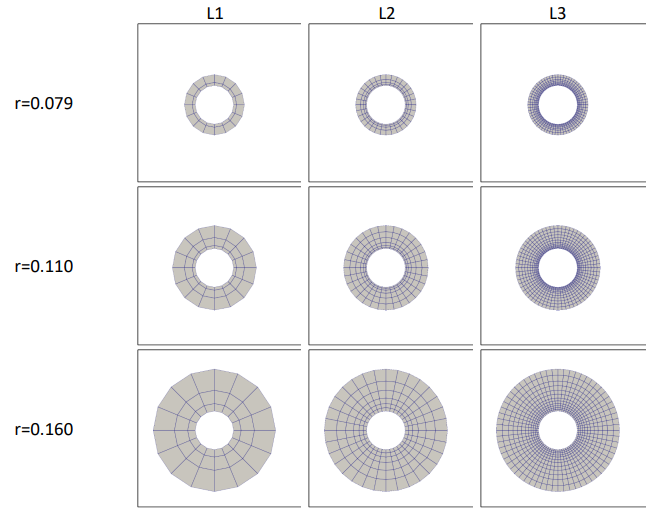} 
        \caption{Body-fitted submeshes on different levels of resolution.}
        \label{fig:sub2}
    \end{subfigure}

    \caption{Background mesh and submeshes for the DFG Benchmark 2D-2.}
    \label{fig:main_figure}
\end{figure}

We performed simulations for the DFG Benchmark 2D-2 using the FBM and
Chimera-S methods. To leverage the advantages of using body-fitted submeshes and gain insights regarding appropriate choices of background mesh/submesh resolutions for our Chimera domain decomposition algorithms, we varied these resolutions and the outer radius of the submesh as shown in Fig. \ref{fig:main_figure}.

\begin{figure}[H]
    \centering
    \includegraphics[width=0.8\textwidth]{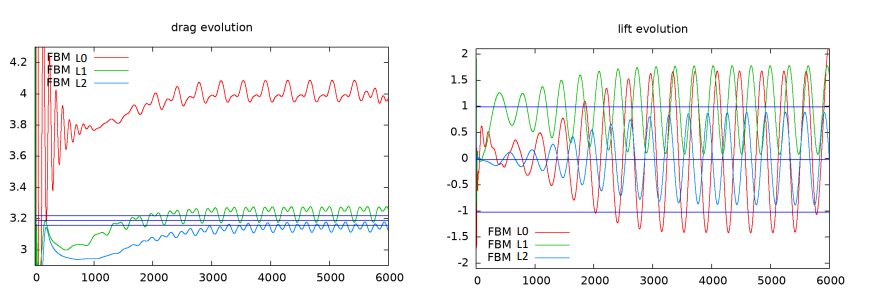} 
    \caption{Evolution of the
      drag and lift coefficients computed with FBM
      on different levels of resolution.}
    \label{fig:drag_standard}
\end{figure}

\begin{figure}[H]
    \centering
    \includegraphics[width=1.0\textwidth]{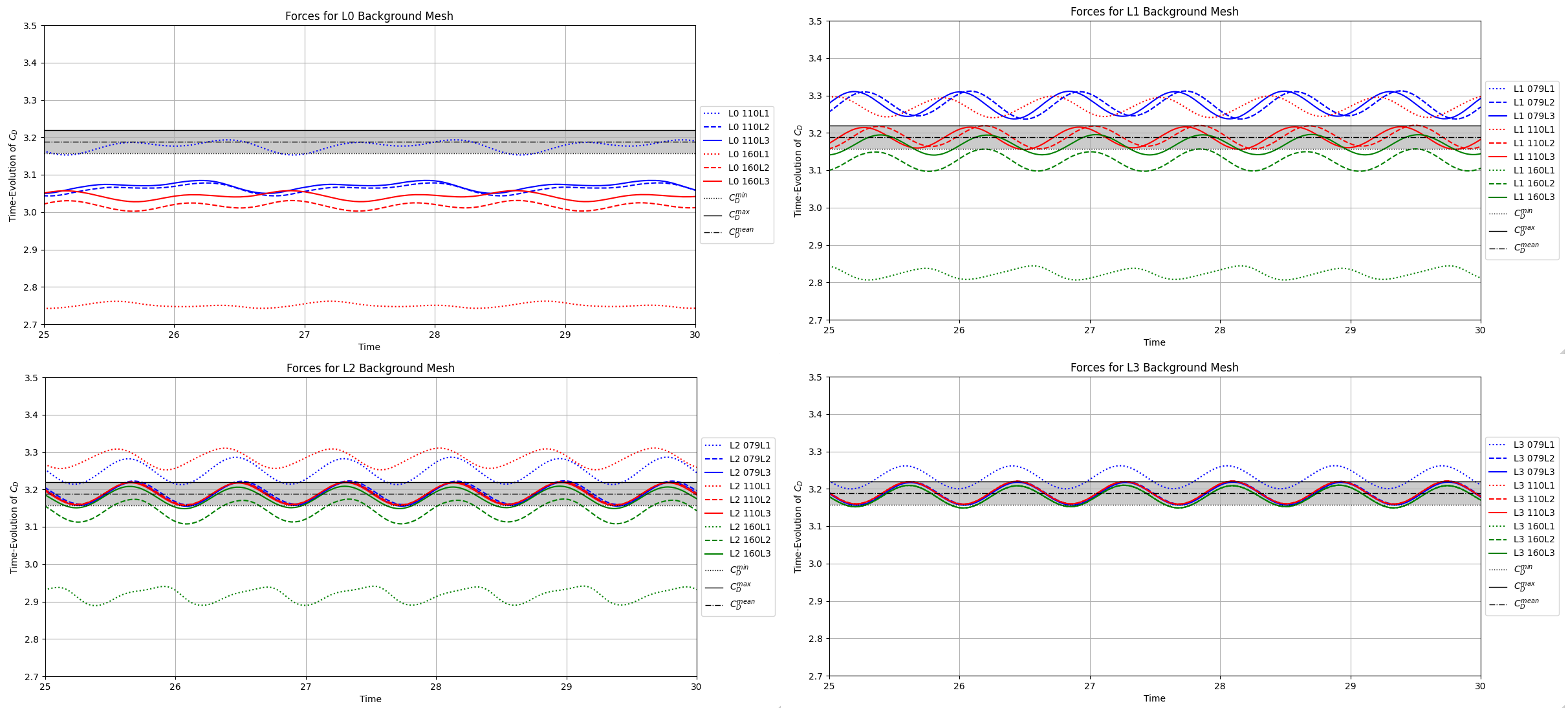} 
    \caption{Evolution of the
      drag coefficient computed with Chimera-S
       on different levels of mesh resolution.}
    \label{fig:drag_chimera}
\end{figure}

\begin{figure}[H]
    \centering
    \includegraphics[width=1.0\textwidth]{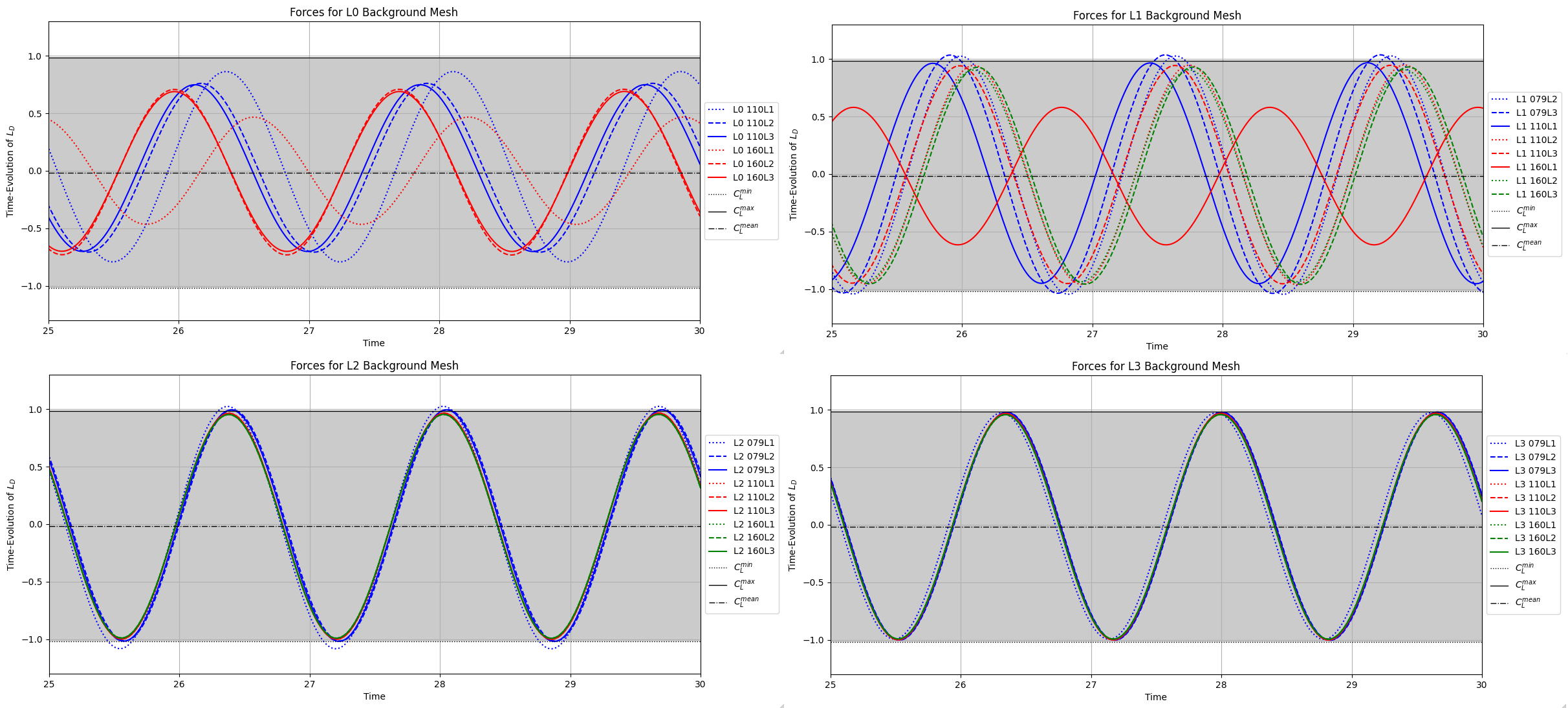} 
    \caption{Evolution of the
      lift coefficient computed with Chimera-S
       on different levels of resolution.}
    \label{fig:lift_chimera}
\end{figure}

The evolution history of the drag and lift coefficients obtained with FBM and Chimera-S is presented in Figure~\ref{fig:drag_standard}-\ref{fig:lift_chimera}. We have marked in the plots the reference values of $C_d^{\min}$, $C_d^{\max}$, $C_d^{\rm{mean}}$, $C_l^{\min}$, $C_l^{\max}$ and $C_l^{\rm{mean}}$, as well as the evolution zone that is defined in previous publications of our research group \cite{dfg,featflow}. The legend for the cases under investigation is as follows: L1\_110L2 corresponds to background mesh resolution L1, submesh radius $r=0.110$, and submesh resolution L2. From the simulation results, we clearly see that all of these three parameters have an influence on the accuracy of the computation. The Chimera method achieves smoother force evolution than FBM and yields more accurate approximations, demonstrating its potential for improved efficiency and reliability.


\subsection{Flow around a moving cylinder}\label{FAC_Moving}

In our second test, a cylinder oscillates horizontally in a
3D channel.
The displacement of its center $(X_c,Y_y,Z_c)$
is sinusoidal in the $x$-direction, while the $y$ and $z$ coordinates remain fixed:
\[
X_c(t) = X_0 + A \sin(2\pi f t), \quad Y_c(t) = Y_0, \quad Z_c(t) = Z_0.
\]
The oscillation amplitude $A$,
oscillation frequency $f$, initial position $(X_0, Y_0, Z_0)$,
and
further parameters of the moving cylinder test are defined in Tables \ref{table1}
and \ref{table2}.

\begin{table}[h!]
\centering
\caption{Cylinder and domain parameters}
\begin{tabular}{ll}
\toprule
Parameter & Value \\
\midrule
Cylinder diameter $D$ & $0.1$ \\
Initial position $(X_0, Y_0, Z_0)$ & $(1.1,\ 0.2,\ 0.1025)$ \\
Oscillation amplitude $A$ & $0.25$ \\
Oscillation frequency $f$ & $0.25$ \\
Domain dimensions & $2.2 \times 0.41 \times 0.1025$ \\
\bottomrule
\end{tabular}
\label{table1}
\end{table}

\begin{table}[h!]
\centering
\caption{Fluid properties}
\begin{tabular}{ll}
\toprule
Property & Value \\
\midrule
Density $\rho$ & $1\, \mathrm{kg/m^3}$ \\
Kinematic viscosity $\nu$ & $10^{-3}\, \mathrm{m^2/s}$ \\
Initial velocity & zero \\
\bottomrule
\end{tabular}
\label{table2}
\end{table}

The primary quantities of interest are again the time-dependent drag
and lift coefficients. Since the constant coordinate $Y_c=0.2$ does not coincide
with the halved height $H=0.205$ of the channel, a non-zero lift force
is exerted on the cylinder. As a reference for the moving cylinder test, we use
the results of 2D simulations conducted by Wan et al~\cite{WanTurek2007a}
on a mesh fitted to the boundary of the cylinder. To obtain the reference
values of the drag and lift coefficients \eqref{eq:cd_cl} for the 3D version,
the results reported in \cite{WanTurek2007a} are scaled by  $\frac{1}{T}$, where $T=0.1025$ is the thickness of the channel.

\begin{figure}[H]
    \centering
    \includegraphics[width=0.9\textwidth]{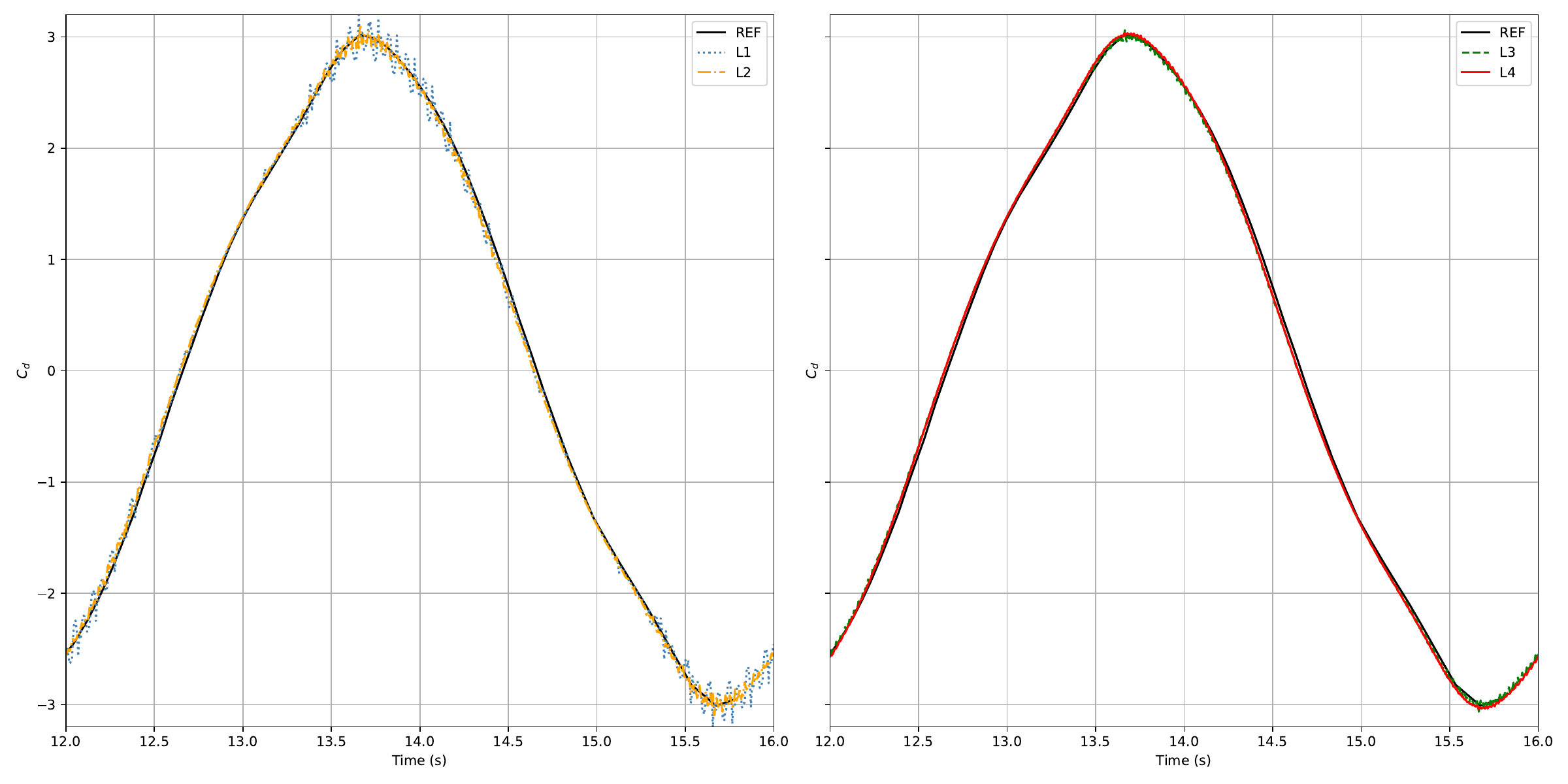} 
    \caption{Evolution of the drag coefficient, Chimera-S method vs. reference data from~\cite{WanTurek2007a}.}
    \label{fig:mc_drag_Codina}
\end{figure}

\begin{figure}[H]
    \centering
    \includegraphics[width=0.9\textwidth]{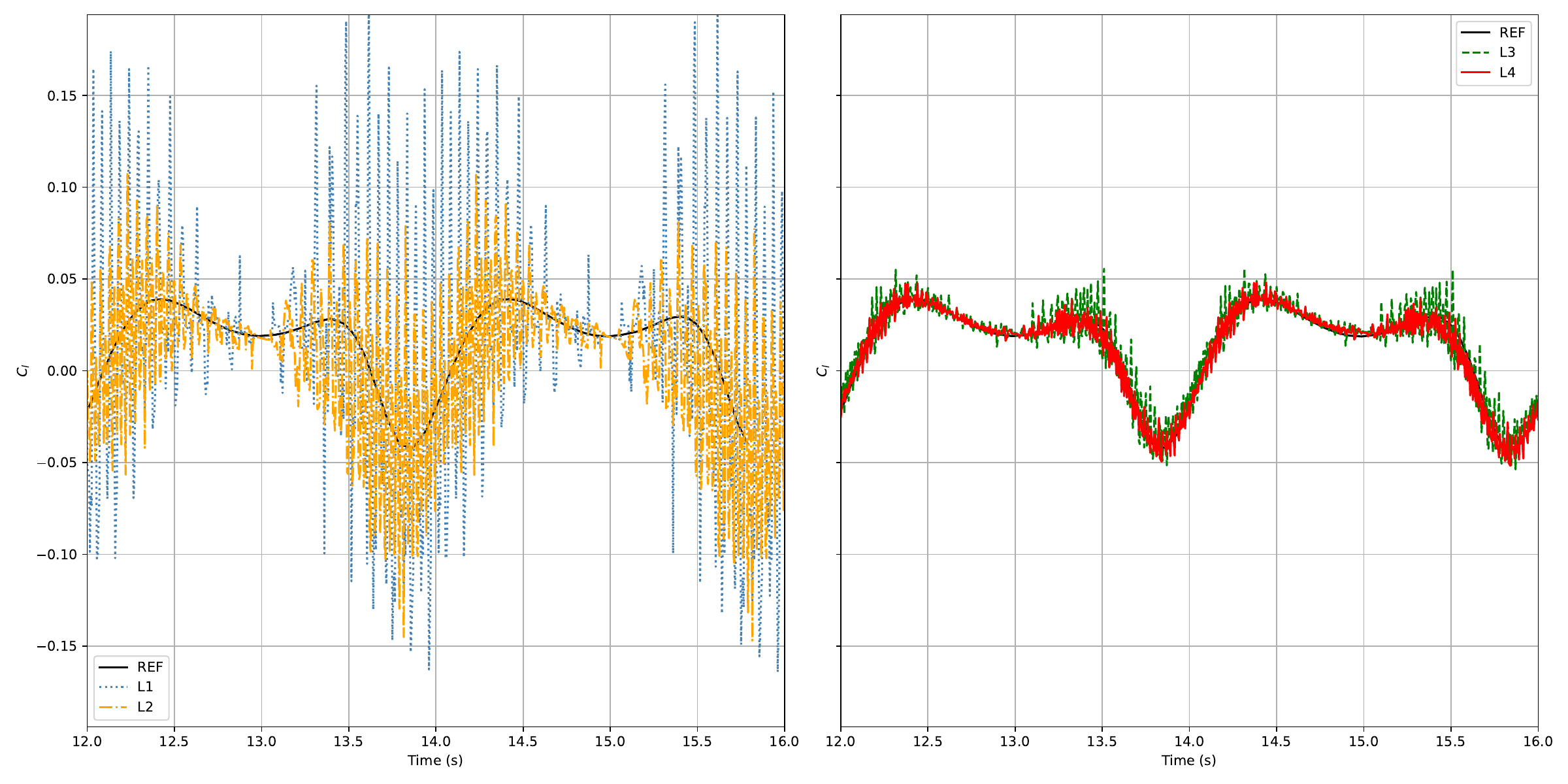} 
    \caption{Evolution of the lift coefficients, Chimera-S method vs. reference data from~\cite{WanTurek2007a}.}
    \label{fig:mc_lift_Codina}
\end{figure}

\begin{figure}[H]
    \centering
    \includegraphics[width=0.9\textwidth]{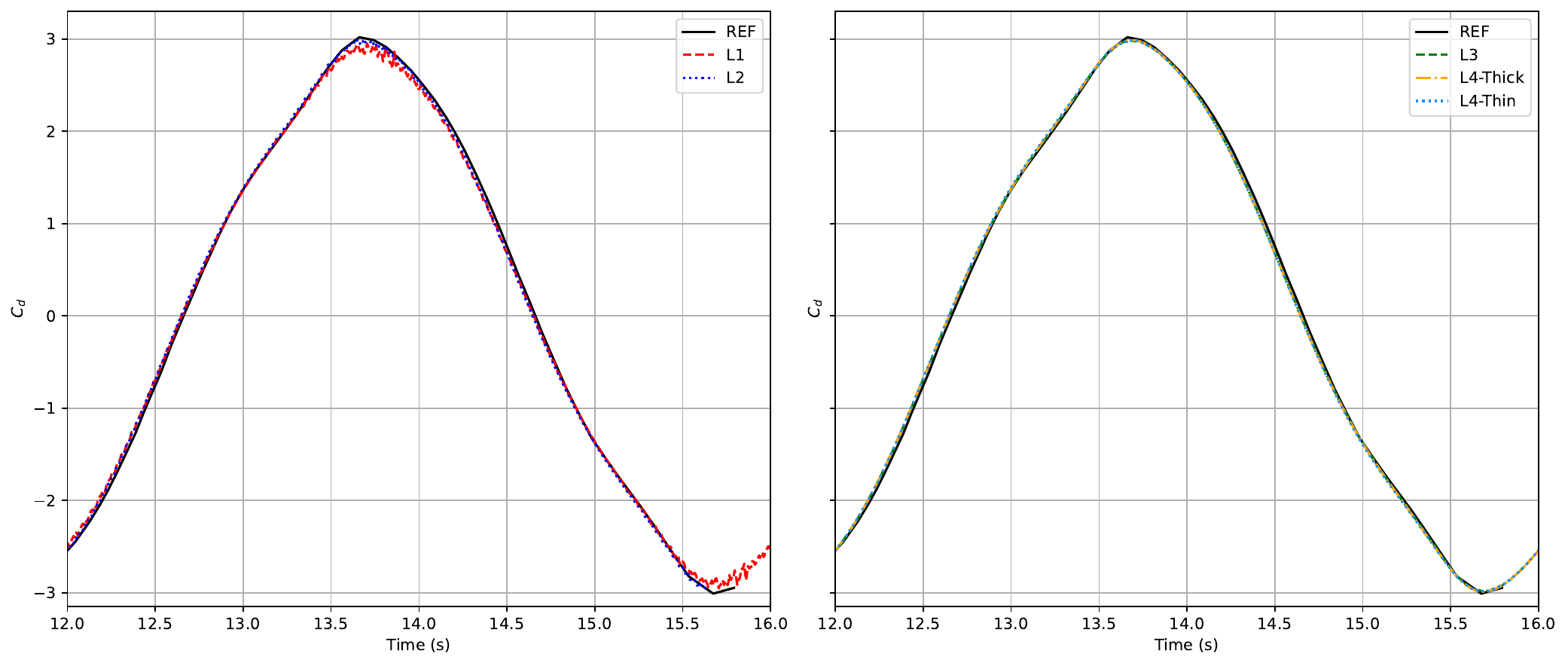} 
    \caption{Evolution of the drag coefficient, Chimera-W method vs. reference data from~\cite{WanTurek2007a}.}
    \label{fig:mc_drag_Penalty}
\end{figure}

\begin{figure}[H]
    \centering
    \includegraphics[width=0.9\textwidth]{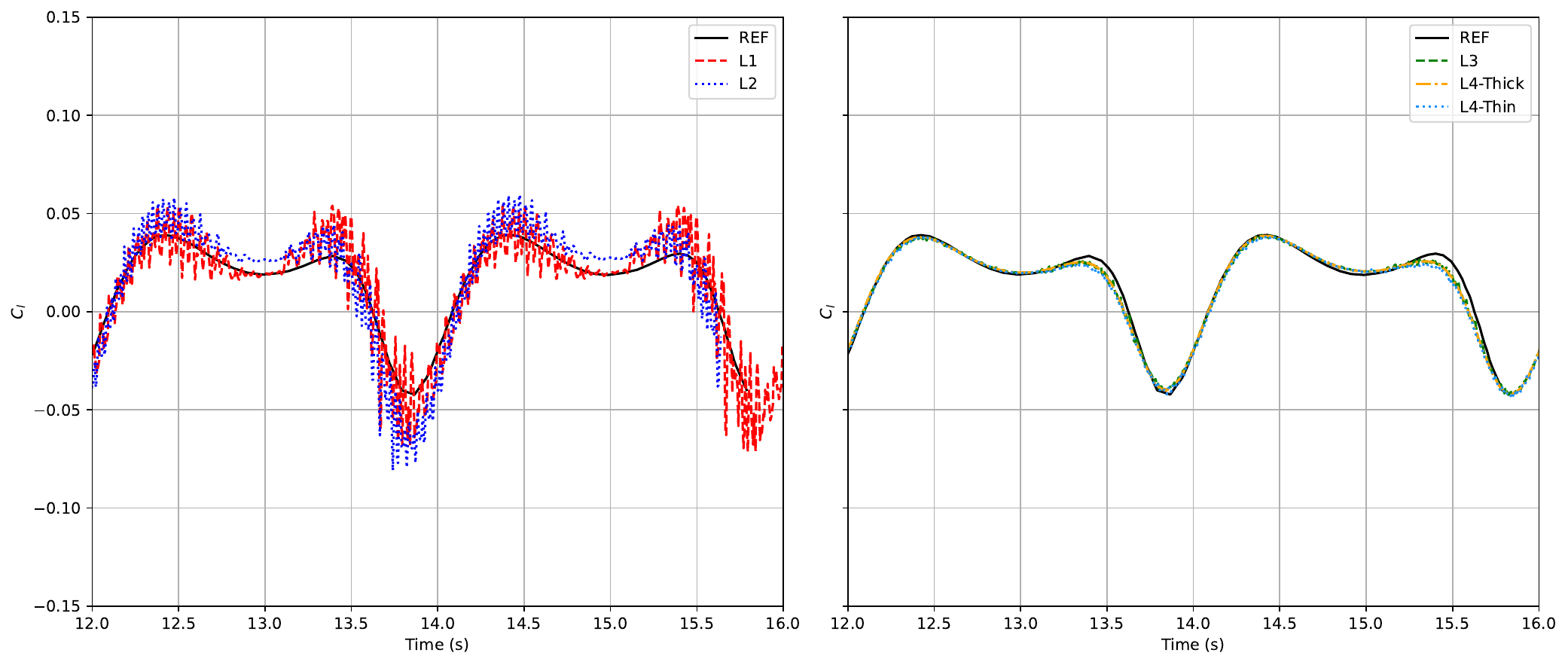} 
    \caption{Evolution of the lift coefficient, Chimera-W method vs. reference data from~\cite{WanTurek2007a}.}
    \label{fig:mc_lift_Penalty}
\end{figure}

We performed mesh refinement studies for both versions (Chimera-S
and Chimera-W) of the domain decomposition method. We plot the
corresponding drag and lift coefficients vs. the reference values.
The Chimera-S results
are shown in Figures \ref{fig:mc_drag_Codina} and \ref{fig:mc_lift_Codina}. They
exhibit strong oscillations because small displacements of the submesh may activate or deactivate a strongly imposed Dirichlet constraint at a node of the background mesh.
The evolution history of the drag and lift coefficients for the Chimera-W
version is displayed in Figures \ref{fig:mc_drag_Penalty} and \ref{fig:mc_lift_Penalty}.
The weak imposition of Dirichlet constraints yields a smoother prediction to the evolving components of the hydrodynamic force. The remaining fluctuations on coarse resolution levels can be attributed to the use of standard numerical quadrature rules for the interior penalty term of the background mesh problem. Similarly to the fringe/hole nodes of the Chimera-S version, small displacements of the submesh may activate or deactivate some quadrature points in an abrupt manner. A possible remedy to this problem is the use of adaptive numerical integration.

Another important finding 
is that even the non-adaptive
implementation of the Chimera-W method exhibits  good
quantitative convergence behavior upon refinement. On the fine-mesh
levels L3 and L4, the results match the reference values within ±1 \% for the drag and ±5 \% for the lift over the entire oscillation cycle. On coarser meshes, Chimera-W can predict the mean drag correctly (because the leading-order contribution is large), but capturing the lift with the same precision
requires mesh resolution that is finer by at least an order of magnitude. This is because the lift is two orders of magnitude less than the drag, so a 1\% error in $C_d$ would correspond to a 100\% error in $C_l$. 

\subsection{Segre--Silberberg migration}
\FloatBarrier
Vertical migration of particles in shear flow is known as the Segre--Silberberg effect. In particular, it is observed in Poiseuille flows that transport particles through a tube or between parallel plates. In this example, particles do not remain near the centerline or the walls. Instead, they migrate to a stable equilibrium position. The typical displacement from the centerline is between $0.2H$ and $0.6H$, where $H$ is the channel radius or half-height. The Segre--Silberberg migration occurs due to a balance between two opposing lift forces: the shear-gradient lift force, which pushes particles away from the center and is caused by velocity gradients, and the wall-induced lift force, which pushes particles away from the walls and is caused by fluid-wall interactions. The stable equilibrium position of a particle significantly depends on the Reynolds number, particle size, and wall-dependent parameters.

\begin{figure}[H]
    \centering
\begin{tikzpicture}

  \draw[thick] (0,0) rectangle (10,4);

  \fill[pattern=crosshatch, pattern color=black!50] (0,3.8) rectangle (10,4);

  \fill[pattern=crosshatch, pattern color=black!50] (0,0) rectangle (10,0.2);
  
  \draw[dashed, thick] (0,2) -- (10,2);
  
  \draw[<->] (8,2) -- (8,4) node[midway, left] {$H$};

  \draw[<->] (0,-0.5) -- (10,-0.5) node[midway, below] {$L$};  
  
  \draw[thick] (0,0.1) .. controls (1.2,2) .. (0,3.9);
  \draw[->, thick, blue] (0.0, 3.0) -- (0.5, 3.0);
  \draw[->, thick, blue] (0.0, 2.0) -- (0.9, 2.0);
  \draw[->, thick, blue] (0.0, 1.0) -- (0.5, 1.0);
  
  \fill[black] (3,3) circle[radius=0.2];

  \draw[<->] (2.7,2.8) -- (2.7,3.2);
  \node[left] at (2.5,3) {$d = 2a$};
  
  \draw[<->] (3.5,2.0) -- (3.5,3.0);
  \node[right] at (3.5,2.5) {$r$};  
  \draw[-] (3.2,3.0) -- (3.5,3.0);  

  \node at (1.8,2.6) {$\bar{a} = \frac{a}{H}$};


\end{tikzpicture}

\caption{Schematic setup of the Segre--Silberberg migration test.}
\label{segre}
\end{figure}
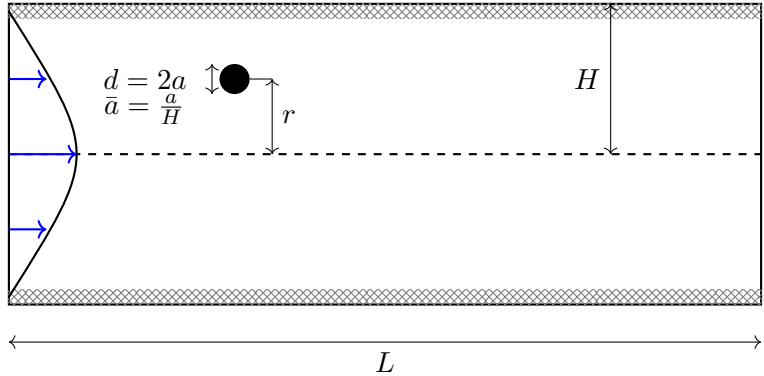

In-depth numerical studies of the Segre--Silberberg effect were performed by Yang et al. \cite{yang2004} for a particle in a tubular channel and a circular particle between parallel planes. In our simulations with the Chimera-S and Chimera-W domain decomposition methods, we compute the equilibrium position of a particle for a Poiseuille flow between two parallel plates. The schematic setup is shown
in Fig.~\ref{segre}. In the corresponding test, the authors of  \cite{yang2004} found that the equilibrium position of the particle is closer to the centerline of the channel than to the nearest wall. This indicates that the wall-induced lift force dominates, pushing the particle away from high shear regions near the wall.

We calculated normalized equilibrium positions $\bar{r}_e=\frac{r_e}{H}$, where $H$ is the half-height of the channel, for two different particle sizes $\bar{a} = \frac{a}{H}$ and a range of Reynolds numbers. Additionally, we studied the dependence on the channel length in our experiments. The results obtained with the Chimera-S and Chimera-W methods are summarized in Tables~\ref{tab:2a=0.10_Codina}-\ref{tab:2a=0.15_Penalty}. The notation for mesh refinement levels follows the same pattern as the legends we used in the description of our numerical results for the DFG cylinder test (see Section \ref{FAC} and, in particular, Fig. \ref{fig:main_figure}). Computations with Chimera-S were performed only for three subsequent resolution levels. Numerical studies of Chimera-W were extended to higher resolution levels (although the submeshes were not refined beyond L3). It can be seen that the two methods are converging to the same results. A comparison with the reference solutions from \cite{yang2004} is provided in Tables \ref{tab:ValidationSegreSilberbergA} and \ref{tab:ValidationSegreSilberbergB}. In order to study the dependence of the equilibrium position on the length $L$ of the computational domain, we ran simulations for $L\in\{5H,10H,20H\}$ on the resolution level~L1. The results of our channel length studies are shown in Table \ref{tab:ChanneLengthSegreSilberberg}. They indicate that the minimum channel length to be used in numerical studies of the Segre--Silberberg effect is $L_{\min}=10H$.

\begin{table}[h]
  \centering
  \vspace{0.75cm}
\begin{tabular}{ccccc}
\toprule
{Re} & \makecell{Background\\mesh\\resolution} & \makecell{Submesh\\resolution} & \makecell{Equilibrium\\position} & \%err \\
\midrule
\multirow{3}{*}{12}
& L1 & 160L1 & 0.40404 & 0.42 \\
& L2 & 130L2 & 0.40654 & 0.19 \\
& L3 & 110L3 & \textbf{0.405}76 & \\
\midrule
\multirow{3}{*}{80}
& L1 & 160L1 & 0.21944 & 1.48 \\
& L2 & 130L2 & 0.21554 & 0.32 \\
& L3 & 110L3 & \textbf{0.216}24 & \\
\midrule
\multirow{3}{*}{180}
& L2 & 130L2 & 0.16669 & 0.65 \\
& L3 & 110L3 & 0.16704 & 0.86 \\
& L4 & 060L3 & \textbf{0.165}62 & \\
\bottomrule
\end{tabular}
\vskip0.5cm

\caption{Equilibrium positions computed with Chimera-S for the
 particle size $2\bar{a}=0.10$.}
\label{tab:2a=0.10_Codina}
\end{table}

\begin{table}[H]
\centering
\begin{tabular}{ccccc}
\toprule
{Re} & \makecell{Background\\mesh\\resolution} & \makecell{Submesh\\resolution} & \makecell{Equilibrium\\position} & \%err \\
\midrule
\multirow{3}{*}{18}
& L1 & 130L1 & 0.44584 & 0.21 \\
& L2 & 110L2 & 0.44676 & 0.00 \\
& L3 & 079L3 & \textbf{0.446}78 & \\
\midrule
\multirow{3}{*}{45}
& L1 & 130L1 & 0.34320 & 0.69 \\
& L2 & 110L2 & 0.34490 & 0.20 \\
& L3 & 079L3 & \textbf{0.345}60 & \\
\midrule
\multirow{3}{*}{180}
& L1 & 130L1 & 0.22008 & 0.07 \\
& L2 & 110L2 & 0.22002 & 0.10 \\
& L3 & 079L3 & \textbf{0.220}24 & \\
\bottomrule
\end{tabular}
\vskip0.5cm

\caption{Equilibrium positions computed with Chimera-S for the
 particle size $2\bar{a}=0.15$.}
\label{tab:2a=0.15_Codina}
\end{table}

\begin{table}[h]
\centering
\begin{tabular}{cccccc}
\toprule
{Re} & \makecell{Background\\mesh\\resolution} & \makecell{Submesh\\resolution} & \makecell{Equilibrium\\position} & \%err & \makecell{\%difference\\(vs. Chimera-S)} \\
\midrule
\multirow{5}{*}{12}
& L1 & 160/L1 & 0.40112 & 0.58 & 1.14 \\
& L2 & 160/L2 & 0.39952 & 0.97 & 1.54 \\
& L3 & 160/L3 & 0.40298 & 0.11 & 0.69 \\
& L4 & 160/L3 & 0.40122 & 0.55 & 1.12 \\
& L4 & 110/L3 & 0.40344 &       & 0.57 \\
\midrule
\multirow{5}{*}{80}
& L1 & 160/L1 & 0.21700 & 0.73 & 0.35 \\
& L2 & 160/L2 & 0.21334 & 0.97 & 1.34 \\
& L3 & 160/L3 & 0.21286 & 1.19 & 1.56 \\
& L4 & 160/L3 & 0.21274 & 1.24 & 1.62 \\
& L4 & 110/L3 & 0.21542 &       & 0.38 \\
\midrule
\multirow{5}{*}{180}
& L1 & 160/L1 & 0.12324 & 26.19 & 25.59 \\
& L2 & 160/L2 & 0.16533 & 0.98 & 0.17 \\
& L3 & 160/L3 & 0.16619 & 0.47 & 0.34 \\
& L4 & 160/L3 & 0.16566 & 0.79 & 0.02 \\
& L4 & 110/L3 & 0.16697 &       & 0.82 \\
\bottomrule
\end{tabular}
\vskip0.25cm

\caption{Equilibrium positions computed with Chimera-S for the
 particle size $2\bar{a}=0.10$.}
\label{tab:2a=0.10_Penalty}
\end{table}


\begin{table}[h]
\centering
\begin{tabular}{cccccc}
\toprule
{Re} & \makecell{Background\\mesh\\resolution} & \makecell{Submesh\\resolution} & \makecell{Equilibrium\\position} & \%err & \makecell{\%difference\\(vs. Chimera-S)} \\
\midrule
\multirow{5}{*}{18}
& L1 & 160/L1 & 0.44252 & 0.88 & 0.95 \\
& L2 & 160/L2 & 0.44206 & 0.98 & 1.06 \\
& L3 & 160/L3 & 0.44300 & 0.77 & 0.85 \\
& L4 & 160/L3 & 0.44318 & 0.73 & 0.81 \\
& L4 & 110/L3 & 0.44644 &       & 0.08 \\
\midrule
\multirow{5}{*}{45}
& L1 & 160/L1 & 0.33550 & 2.36 & 2.92 \\
& L2 & 160/L2 & 0.33848 & 1.49 & 2.06 \\
& L3 & 160/L3 & 0.33920 & 1.28 & 1.85 \\
& L4 & 160/L3 & 0.33904 & 1.33 & 1.90 \\
& L4 & 110/L3 & 0.34360 &       & 0.58 \\
\midrule
\multirow{5}{*}{180}
& L1 & 160/L1 & 0.21756 & 0.44 & 1.22 \\
& L2 & 160/L2 & 0.21600 & 1.15 & 1.93 \\
& L3 & 160/L3 & 0.21544 & 1.41 & 2.18 \\
& L4 & 160/L3 & 0.21536 & 1.45 & 2.22 \\
& L4 & 110/L3 & 0.21852 &       & 0.78 \\
\bottomrule
\end{tabular}
\vskip0.5cm

\caption{Equilibrium positions computed with Chimera-W for the
 particle size $2\bar{a}=0.15$.}
\label{tab:2a=0.15_Penalty}
\end{table}

\begin{table}[h]
\centering
\begin{tabular}{c c c c c c c}
\toprule
$2\bar{a}$ & $Re$ & Short & Medium & Long &
{M.\,vs.\,S [\%]} & {M.\,vs.\,L [\%]}\\
\midrule
\multirow{3}{*}{0.10}
 & 18  & 0.3330 & 0.3427 & 0.3426 &  -2.9 &  -0.1 \\
 & 80  & 0.2175 & 0.2194 & 0.2189 &  -0.9 &  -0.2 \\
 & 180 &-0.1061 & 0.1667 & 0.1579 & -163.6 &  -5.3 \\
\midrule
\multirow{3}{*}{0.15}
 & 18  & 0.4437 & 0.4458 & 0.4458 &  -0.5 &   0.0 \\
 & 45  & 0.3349 & 0.3432 & 0.3433 &  -2.4 &   0.0 \\
 & 180 & 0.1426 & 0.2201 & 0.2184 & -35.2 &  -0.8 \\
\bottomrule
\end{tabular}
\vskip0.5cm

\caption{Equilibrium positions computed with Chimera-S in
  short ($L=5H$), medium ($L=10H$), and long ($L=20H$) channel simulations
  using the L1 background mesh resolution.}
\label{tab:ChanneLengthSegreSilberberg}

\vskip0.5cm

\centering
\begin{subtable}[t]{0.48\textwidth}
\centering
\begin{tabular}{lcccc}
\toprule
Re & ref. \cite{yang2004} & FBM & Chimera-W & Chimera-S \\
\midrule
12   & 0.413 & 0.388 & 0.403 & {\bf 0.405}8 \\
80   & 0.222 & 0.225 & 0.215 & {\bf 0.216}2 \\
180  & 0.174 & 0.174 & 0.167 & {\bf 0.165}6 \\
\bottomrule
\end{tabular}
\caption{Particle size $2\bar{a} = 0.10$.}
\label{tab:ValidationSegreSilberbergA}
\end{subtable}
\vskip0.5cm

\begin{subtable}[t]{0.48\textwidth}
\centering
\begin{tabular}{lcccc}
\toprule
Re & ref. \cite{yang2004} & FBM & Chimera-W & Chimera-S \\
\midrule
18   & 0.454 & 0.450 & 0.446 & {\bf 0.446}8 \\
45   & 0.359 & 0.351 & 0.344 & {\bf 0.345}6 \\
180  & 0.234 & 0.227 & 0.219 & {\bf 0.220}2 \\
\bottomrule
\end{tabular}
\caption{Particle size $2\bar{a} = 0.15$.}
\label{tab:ValidationSegreSilberbergB}
\end{subtable}
\vskip0.5cm

\caption{Converged equilibrium positions, Chimera-W and
  Chimera-S vs. reference
  data from \cite{yang2004}.}
\label{tab:ValidationSegreSilberbergCombined}
\end{table}

\begin{figure}[h]
  \centering
  \begin{subfigure}[t]{0.3\textwidth}
    \centering
    \includegraphics[width=\linewidth]{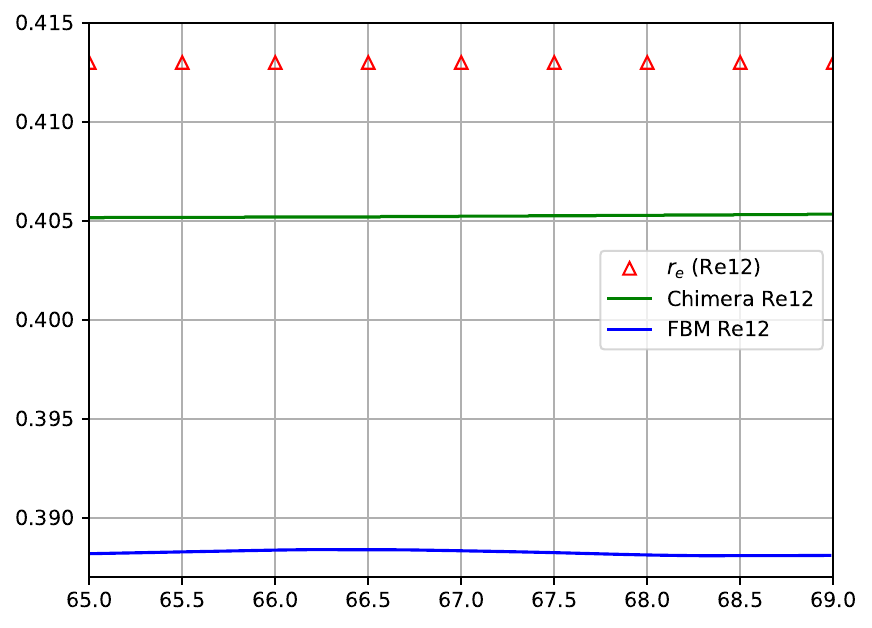}
    \caption{Re=12}
    \label{fig:plot1}
  \end{subfigure}
  \hfill
  \begin{subfigure}[t]{0.3\textwidth}
    \centering
    \includegraphics[width=\linewidth]{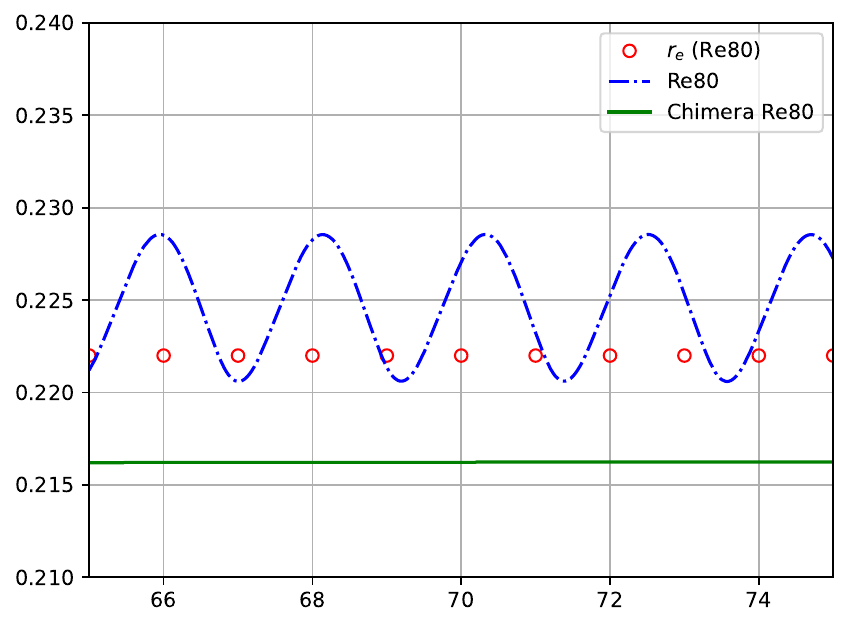}
    \caption{Re=80}
    \label{fig:plot2}
  \end{subfigure}
  \hfill
  \begin{subfigure}[t]{0.3\textwidth}
    \centering
    \includegraphics[width=\linewidth]{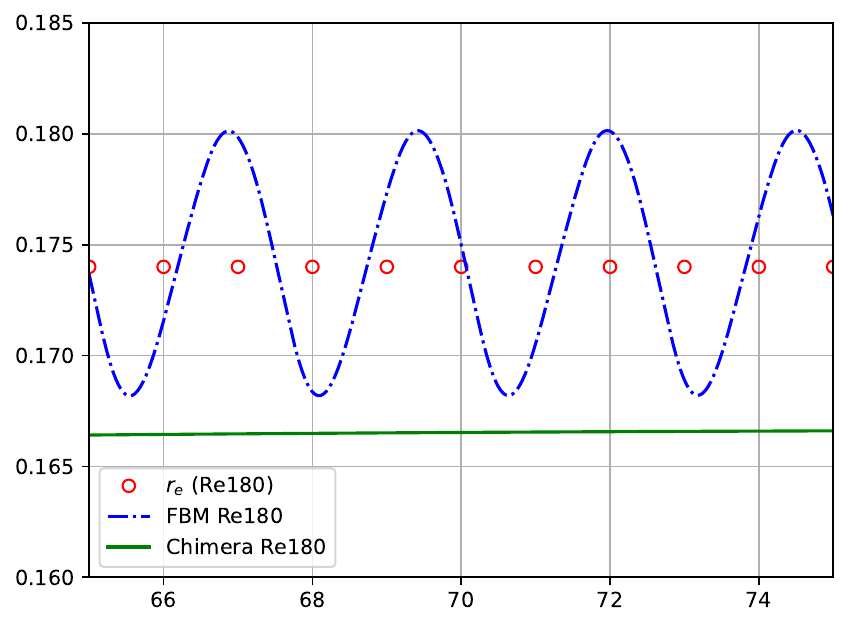}
    \caption{Re=180}
    \label{fig:plot3}
  \end{subfigure}

  \vspace{1em} 

  \begin{subfigure}[t]{0.3\textwidth}
    \centering
    \includegraphics[width=\linewidth]{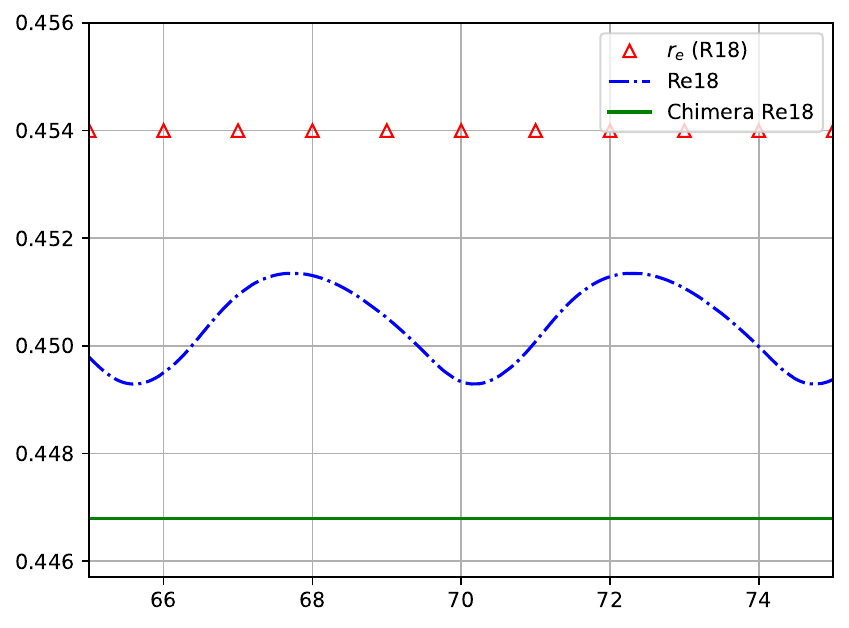}
    \caption{Re=18}
    \label{fig:plot4}
  \end{subfigure}
  \hfill
  \begin{subfigure}[t]{0.3\textwidth}
    \centering
    \includegraphics[width=\linewidth]{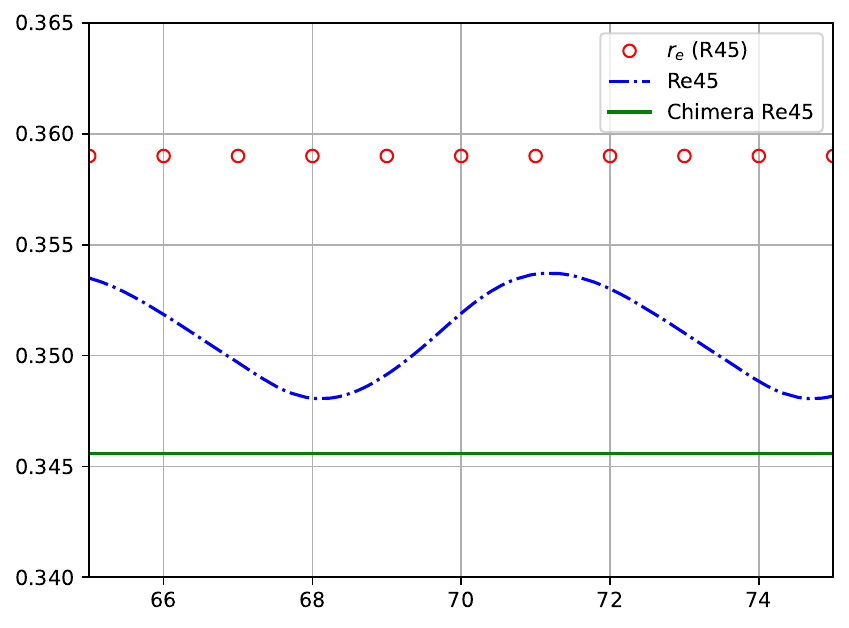}
    \caption{Re=45}
    \label{fig:plot5}
  \end{subfigure}
  \hfill
  \begin{subfigure}[t]{0.3\textwidth}
    \centering
    \includegraphics[width=\linewidth]{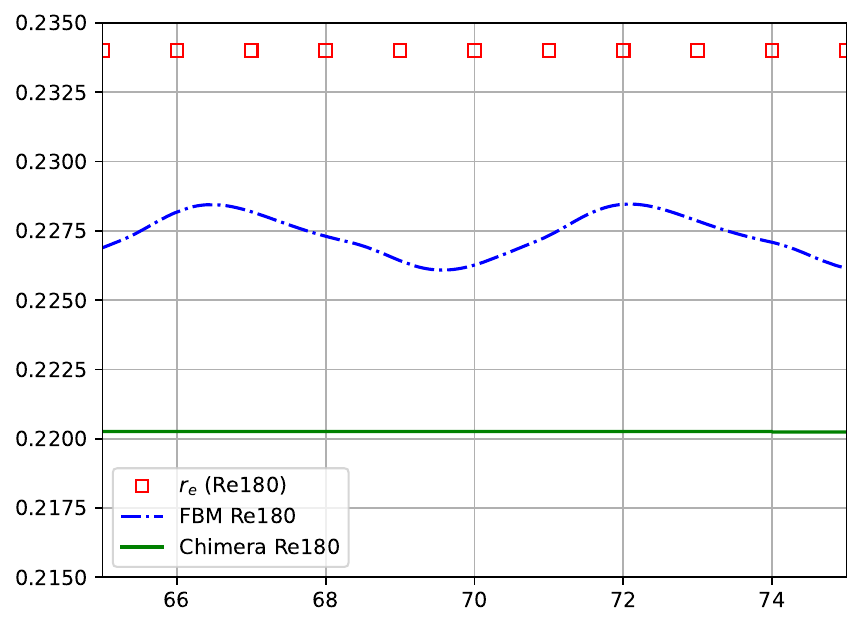}
    \caption{Re=180}
    \label{fig:plot6}
  \end{subfigure}
  \vskip0.5cm
  
  \caption{Stability and oscillations of Segre–Silberberg equilibrium positions,
    FBM vs. Chimera-S.}
  \label{fig:all_plots}
\end{figure}

To gain further insights into the numerical behavior of Chimera submesh methods,
we simulated the Segre--Silberberg migration effect using the classical single-mesh implementation of FBM. A visual comparison of the results presented in Fig. \ref{fig:all_plots} indicates that the Chimera-S solutions are significantly less oscillatory in five of the six cases under investigation. The FBM version converges to a value that is close to the reference solution but exhibits a periodic oscillation of $O(10^{-3})$ to $O(10^{-2})$ around its mean value. The amplitude of the oscillation seems to be correlated with the Reynolds number.

The Chimera-W results (not presented here) are not quite as accurate as the Chimera-S results in this particular test. The superior performance of Chimera-S can be attributed to the fact that the set of fringe and hole nodes reaches a steady state in the equilibrium position. The performance of Chimera-W could be further improved by using adaptive quadrature for the interior penalty term.

\section{Conclusions}
\label{sec:concl}

In this work, we extended the Chimera method presented in \cite{codina} to
time-dependent particulate flows and proposed a new implementation of the
Dirichlet--Robin coupling. The results of our numerical experiments indicate
that weak imposition of Dirichlet constraints via interior penalty terms
stabilizes the hydrodynamic forces in simulations with moving particles. The benefits of the proposed Chimera multimesh method include simplicity, efficiency, and continuous dependence of numerical solutions on the location of moving particles. The use of body-fitted submeshes eliminates the need for mesh deformation techniques designed to enhance the accuracy of fictitious boundary / subspace projection algorithms with strongly imposed Dirichlet constraints \cite{steffen,fbm2012,WanTurek2007b}. Extension to overlapping submeshes are feasible and can be performed by adapting the general framework developed in \cite{multimesh2020,multimesh2019} for projection schemes using discontinuous Galerkin weak forms and a different kind of interior penalization.

  \section*{Acknowledgments}
  This work was supported by the German Research Foundation (DFG) under grant
  KU 1530/28-1 (TU 102/77-1). The authors gratefully acknowledge collaboration
   on this project
  with Prof. Yuliya Gorb (National Science Foundation) and Prof. Alexey Novikov
  (Pennsylvania State University).

\bibliographystyle{plain} 
\bibliography{references} 

\section*{Appendix A. Matrices and vectors of the discrete problem}

The finite element approximations to the background velocity and pressure
are given by
$$
\bu_h(\bx,t)=\sum_{j=1}^{N_h} u_j(t)\boldsymbol{\varphi}_j(\bx),
\quad p_h(\bx,t)=\sum_{k=1}^{M_h} p_k(t)\psi_k(\bx),
$$
where $\boldsymbol{\varphi}_j$ and $\psi_k$ are basis functions spanning the
spaces $\mathbf V_h$ and $Q_h$, respectively. In our description of discrete
problems, $u=(u_j)$ and $p=(p_k)$ are vectors containing the
coefficients of the above FE approximations. Recall that the
systems of equations considered
in Section \ref{sec:proj} depend on
$$
A=(a_{ij}),\qquad B=(b_{ik}),\qquad M_C=(m_{ij}),\qquad M_L=(\tilde m_{ij}),
\qquad D=(d_{ij}).
$$
In view of \eqref{intpen} and \eqref{wf}, the entries of these matrices
and the components of $f=(f_i)$ are defined by
\begin{align*}
  a_{ij}&=\frac{m_{ij}}{\Delta t}+\theta\left[
  \rho_f\int_{\Omega}
  (\bu_h\cdot\nabla\boldsymbol{\varphi}_j)\cdot\boldsymbol{\varphi}_i\dx
  +
  \frac{\mu_f}2\int_\Omega \mathbf D(\boldsymbol{\varphi}_j):\mathbf D(\boldsymbol{\varphi}_i)\dx\right],\\
  b_{ik}&=-\int_\Omega \psi_k\nabla\cdot\boldsymbol{\varphi}_i\dx,\qquad
  m_{ij}=\rho_f\int_{\Omega}\boldsymbol{\varphi}_j\cdot\boldsymbol{\varphi}_i\dx,
  \qquad \tilde m_{ij}=\rho_f\boldsymbol{\varphi}_j(\bx_i)
  \cdot \int_{\Omega}\boldsymbol{\varphi}_i\dx,\\
  d_{ij}&=\gamma_{\max}\sum_{k=1}^{N_p}\left[\int_{\hat\Omega_{k,h}}
  \beta_k\boldsymbol{\varphi}_j\cdot\boldsymbol{\varphi}_i\dx+
  \int_{B_{k,h}}\boldsymbol{\varphi}_j\cdot\boldsymbol{\varphi}_i\dx
  \right],
  \\
  f_i&=\sum_{j=1}^{N_h}\left(
\frac{m_{ij}}{\Delta t}-(1-\theta)\left[
  \rho_f\int_{\Omega}
  (\bu_h\cdot\nabla\boldsymbol{\varphi}_j)\cdot\boldsymbol{\varphi}_i\dx
  +
  \frac{\mu_f}2\int_\Omega \mathbf D(\boldsymbol{\varphi}_j):\mathbf D(\boldsymbol{\varphi}_i)\dx\right]
    \right)u_j.
\end{align*}
In the formula for $\tilde m_{ij}$, we denote by
$\mathbf{x}_i$ the nodal point associated with the basis
function $\boldsymbol{\varphi}_i$.

The data of the weakly imposed Dirichlet constraints is built into 
$g=(g_i)$ with
$$
g_i=\gamma_{\max}\sum_{k=1}^{N_p}\left[\int_{\hat\Omega_{k,h}}
  \beta_k\hbu_h\cdot\boldsymbol{\varphi}_i\dx+
  \int_{B_{k,h}}\bU_h\cdot\boldsymbol{\varphi}_i\dx
  \right].
$$

\end{document}